\documentclass{amsart}

\usepackage{amsthm,amssymb,latexsym,amsmath}
\usepackage{mathrsfs}
\usepackage{graphicx}
\usepackage{color}
\usepackage[latin1]{inputenc}

\newcommand{\CC}{{\mathbb C}}

\newcommand{\CIPT}{\mathbb{P}^n}
\newcommand{\tf}{T_{{\mathscr F}}}

\newcommand{\ddd}{\mathcal D}

\renewcommand{\dim}{\mathrm{dim}}

\newcommand{\OO}{\mathcal O}
\newcommand{\F}{\mathscr F}
\newcommand{\fol}{\mathscr F}

\newcommand{\NN}{{\mathbb{N}}}

\newcommand{\PP}{\mathbb{P}}

\newcommand{\singf}{{\rm Sing}({\mathscr F})}

\newcommand{\ZZ}{\mathbb Z}

\newcommand{\Sing}{{\rm Sing}}
\newcommand{\cod}{{\rm cod}}

\newtheorem{lema}{Lemma}[section]

\newtheorem{teo}[lema]{Theorem}
\newtheorem*{teo1}{Theorem 1}
\newtheorem*{teo2}{Theorem 2}
\newtheorem*{teo3}{Theorem 3}

\newtheorem*{cor1}{Corollary 1}

\newtheorem*{teorema}{Theorem}

\newtheorem{prop}[lema]{Proposition}

\theoremstyle{definition}
\newtheorem{remark}[lema]{Remark}
\newtheorem{defi}[lema]{Definition}
\newtheorem{exe}[lema]{Example}

\newenvironment{demostracion}
{\noindent {{\it  Proof.}}}

\begin{document}


\title{Residue Formulas for logarithmic Foliations and applications}


\begin{abstract}
In this work we  prove a  Baum-Bott type formula for non-compact complex manifold of the form $\tilde{X}=X-\ddd$, where $X$ is a complex compact manifold and $\ddd$ is a normal crossing  divisor  on  $X$. As applications, 
we provide a Poincar\'e-Hopf type Theorem and   an optimal  description  for a smooth hypersurface $\ddd$ invariant by an one-dimensional foliation $\fol$ on  $\mathbb{P}^n$ satisfying $\Sing(\fol) \subsetneq \ddd.$
\end{abstract}

\author{Maur\'icio Corr\^ea}

\address{\noindent  Maur\'icio Corr\^ea\\
Departamento de Matem\'atica \\
Universidade Federal de Minas Gerais\\
Av. Ant\^onio Carlos 6627 \\
30123-970 Belo Horizonte -MG, Brazil} \email{mauriciojr@ufmg.br}

\author{Diogo da Silva Machado}

\address{\noindent  Diogo da Silva Machado\\
Departamento de Matem\'atica \\
Universidade Federal de Vi\c cosa\\
Avenida Peter Henry Rolfs, s/n - Campus Universitário \\
36570-900 Vi\c cosa- MG, Brazil} \email{diogo.machado@ufv.br}

\subjclass{Primary 32S65, 32S25, 14C17 } 

\keywords{Logarithmic foliations, Poincar\'e-Hopf type Theorem, residues}


\maketitle



\section{Introduction}

In \cite{BB1} P. Baum and  R. Bott developed a  work about  residues of singularities of  holomorphic foliations on complex manifolds. In the case of  one-dimensional holomorphic foliation $\fol$, with isolate singularities, on an $n$-dimensional complex compact manifold $X$ we have the following classical Baum-Bott formula:
\begin{eqnarray}\label{bbformula}
  \displaystyle \hspace*{2.0 cm}\int_{X}c_n(T_X - \tf) \,\,\, =\,\,\, \sum_{p\in \Sing(\fol)} \mu_{p}(\fol),\,\,\,\,\,\,\,\,\,\,\,\,\,\,\,\,\,\,\,\, {\small{\mbox{({\it Baum-Bott formula}}})}
\end{eqnarray}
\noindent where the $\mu_{p}(\fol)$ are the Milnor number of $\fol$ in $p$.  Baum-Bott formula is a generalization (for holomorphic vector fields)  of the Poincar\'e-Hopf Theorem
$$
\int_{X}c_n(T_X ) \,\,\, =\,\,\, \sum_{p\in \Sing(\fol)} \mu_{p}(\fol),
$$ 
where  $\fol$  is a foliation induced by a global  holomorphic vector field, with isolated singularities,  on  $X$. 

In this work we provide a  Baum-Bott type formula for non-compact complex manifold of the form  $\tilde{X}=X-\ddd$, where $X$ is a complex compact manifold and $\ddd$ is an analytic divisor contained in $X$ invariant by an one-dimensional holomorphic foliation $\fol$ which is called by \textit{logarithmic  foliation} along $\ddd$. 
As an application, we obtained a  Poincar\'e-Hopf type Theorem for these non-compact manifolds. Furthermore, for logarithmic foliations on projective spaces, we prove a necessary and sufficient conditions for all singularities of the foliation occur in an analytic  invariant hypersurface.

We prove the following result.
\begin{teo1}\label{teo1_1} Let $\tilde{X}$ be an $n$-dimensional complex manifold such that $\tilde{X} = X - \ddd$, where $X$ is an $n$-dimensional complex compact manifold and  $\ddd$ is a smooth hypersurface on $X$. Let $\fol$ be a foliation of dimension one on $X$, with isolated singularities and logarithmic along $\ddd$. 
Suppose  that 
$Ind_{log\, \ddd,p}(\fol)=0$,  for all $p \in \Sing(\fol)\cap \ddd$.
Then
$$
\int_{X}c_{n}(T_{X}(-\log\, \ddd)- T_{\fol}) = \sum_{p\in \Sing(\fol)\cap (X\setminus \ddd)}\mu_p(\fol).
$$
\end{teo1}
Here, $Ind_{log\, \ddd,p}(\fol)$ denotes the logarithmic index of $\fol$ on $p$, see  section \ref{seclog}.

 The classical Gauss-Bonnet theorem for a complex  compact manifold  $X$,   proved by S. Chern in \cite{chern},  says  us that 
\begin{equation}\label{class_Gauss}
\displaystyle \int_{X}c_n(T_X)  = \chi (X).
\end{equation}
The following version of Gauss-Bonnet formula for non-compact manifolds was initially proposed by S. Iitaka \cite{Lita} and proved by Y. Norimatsu \cite{YN}, R. Silvotti \cite{RS} and P. Aluffi \cite{Aluffi}:

\begin{teorema}[Norimatsu-Silvotti -Aluffi]
Let $\tilde{X}$ be an $n$-dimensional complex manifold such that $\tilde{X} = X - \ddd$, where $X$ is an $n$-dimensional complex compact manifold and $\ddd$ is a normally crossing hypersurface on $X$. Then 
$$
\int_{X}c_n(T_X(-\log\, \ddd)) \,\,\, =\,\,\, \chi (\tilde{X}),
$$
 where  $\chi (\tilde{X})$ denotes the Euler characteristic given by 
$$
\chi (\tilde{X}) = \displaystyle \sum_{i=1}^n \dim \  H^i_{c}(\tilde{X}, \CC).
$$
\end{teorema}
X. Liao in \cite{Liao} has provided more general  formulas in terms of  Chern-Schwartz-MacPherson class of $\tilde{X}$.

In the Section \ref{sec05}, we consider the case where $\ddd$ is a normal  crossing hypersurface and we prove the following Baum-Bott type formula:

\begin{teo2} \label{teo2_2} Let $\tilde{X}$ be an $n$-dimensional complex manifold such that $\tilde{X} = X - \ddd$, where $X$ is a $n$-dimensional complex compact manifold,  $\ddd$ is a normally crossing hypersurface on $X$. Let $\fol$ be a foliation on $X$ of dimension one, with isolated singularities (non-degenerates) and logarithmic along $\ddd$. Then,
\begin{equation}\label{pp1}
\displaystyle\int_{X}c_{n}(T_{X}(-\log\, \ddd)- T_{\fol}) = \sum_{p\in \Sing(\fol) \cap (\tilde{X})}\mu_p(\fol).
\end{equation}
\end{teo2}

As a  consequence of Theorem $2$ and Norimatsu-Silvotti -Aluffi Theorem  we obtain  the following Poincar\'e-Hopf type Theorem.

\begin{cor1} \label{cor1_1}
Let $\tilde{X}$ be an $n$-dimensional complex manifold such that $\tilde{X} = X - \ddd$, where $X$ is an $n$-dimensional complex compact  manifold, $\ddd$ is a reduced  normal crossing hypersurface on $X$.  Let $\fol$ be a foliation on $X$ of dimension one given by a global holomorphic vector field, with isolated singularities (non-degenerates) and logarithmic along $\ddd$. Then 
$$
 \chi(\tilde{X}) =\sum_{p\in \Sing(\fol)\cap  \tilde{X} }\mu_p(\fol),
$$
where $\mu_p(\fol)$ denotes the Milnor number of $\fol$ on $p$.
\end{cor1}

Finally, in  the Section \ref{sec06}, we prove   a complete characterization in order that an  invariant hypersurface contains all the singularities of the projective foliation.

\begin{teo3} \label{teo_principal}  Let $\ddd \subset \CIPT$ be a smooth and irreducible hypersurface and let $\fol$ be a foliation of dimension one on $\CIPT$ , with isolated singularities (non-degenerates) and logarithmic along $\ddd$. Then, the following proprieties holds
\begin{itemize}

\item [{\bf (1)}] If $n$ is odd, then:\\

\begin{itemize}

\item [{\bf (a)}] $\# \left [\singf \cap {\CIPT \backslash \ddd}\right ]>0 \,\,\,\Longleftrightarrow\,\,\, \deg(\ddd)<\deg(\fol)+1$;
\\
\item [{\bf (b)}] $\# \left [\singf \cap {\CIPT \backslash \ddd}\right ]=0\,\,\, \Longleftrightarrow \,\,\,\deg(\ddd)=\deg(\fol)+1$.
\\
\end{itemize}

\item [{\bf (2)}] If $n$ is even, then:
\begin{itemize}
\item [{\bf (a)}] $\# \left [\singf \cap {\CIPT \backslash \ddd}\right ] >0 \,\,\,\Longleftrightarrow\,\,\,\left \{ \begin{array}{l}
  \deg(\ddd)\neq \deg(\fol)+1\\
\mbox{or}\\
\deg(\ddd)=\deg(\fol)+1,\,\,\mbox{with}\,\, \deg(\fol)\neq 0 \\
\end{array}\right.$
\\
\item [{\bf (b)}] $\# \left [\singf \cap {\CIPT \backslash \ddd}\right ] = 0 \,\,\,\Longleftrightarrow\,\,\, \deg(\ddd) = 1 \,\,\mbox{and}\,\, \deg(\fol)= 0.$
\\
\end{itemize}

\item[{\bf (3)}] In general, we have the formula 
\begin{eqnarray}\nonumber
\# \left [\singf \cap {\CIPT \backslash \ddd}\right ] = \displaystyle\sum^{n}_{i=0}(-1)^i(\deg(\ddd)  -1)^i\deg(\fol)^{n-i}.
\end{eqnarray}

\end{itemize}

\end{teo3}

Observe that if   $n$ is odd, then $\Sing(\fol) \subsetneq \ddd$ if and only if the Soares's bound for  the Poincar\'e problem  is  achieved \cite{MGS}.

\subsection*{Acknowledgments}
We are grateful to Gilcione Nonato, Jean-Paul Brasselet,  Tatsuo Suwa and  Marcio G.  Soares    for interesting conversations.
This work was partially supported by CNPq, CAPES, FAPEMIG and FAPESP-2015/20841-5. Finally, we would like to thank the referees by the suggestions, comments and improvements to the exposition.

\section{Preliminaries}

\subsection{Logarithmics forms and logarithmics vector fields}

\hyphenation{theo-ry}

Let $X$ an $n$-dimensional complex manifold and $\ddd$ a reduced hypersurface on $X$. Given a  meromorphic $q$-form $\omega$ on $X$, we say that $\omega$ is a {\it logarithmic} $q$-form along $\ddd$ at $x\in X$ if the following conditions occurs:

\begin{enumerate}
\item[(i)] $\omega$ is holomorphic on $X - \ddd$;
\medskip
\item[(ii)] If $h=0$ is a reduced equation of $\ddd$, locally at $x$, then $h \,\omega$ and $h\, d \omega$ are holomorphic.
\end{enumerate}
\noindent Denoting by $\Omega_{X,x}^q(\log\, \ddd)$ the set of germs of logarithmic $q$-form  along $\ddd$ at $x$, we define the following coherent sheaf of $\OO_X$-modules
$$
\Omega_{X}^q(\log\, \ddd):= \bigcup_{x\in X} \Omega_{X,x}^q(\log\, \ddd),
$$
which is called  by {\it sheaf of logarithmics $q$-forms along $\ddd$}.  See \cite{Deligne}, \cite{Katz} and  \cite{Sai}  for details.

Now, given $x\in X$, let $v\in T_{X,x}$ be germ at $x$ of a holomorphic vector field on $X$. We say that $v$ is a {\it logarithmic vector field along of $\ddd$ at $x$}, if $v$ satisfies the following condition: if $h=0$ is a equation of $\ddd$, locally at $x$, then the derivation $v(h)$ belongs to the ideal $\langle h_x\rangle\OO_{X,x}$. Denoting by $T_{X,x}(-\log\, \ddd)$ the set of germs of logarithmic vector field along of $\ddd$ at $x$, we define the following coherent sheaf of $\OO_X$-modules
$$
T_{X}(-\log\, \ddd):= \bigcup_{x\in X}T_{X,x}(-\log\, \ddd),
$$
which is called by  {\it sheaf of logarithmics vector fields along $\ddd$}. 

 It is known that $\Omega_{X}^1(\log\, \ddd)$ and $T_{X}(-\log \,\ddd)$ is always a reflexive sheaf, see  \cite{Sai} for more details. If  $\ddd$ is an analytic hypersurface with normal crossing singularities, the sheaves $\Omega_{X}^1(\log\, \ddd)$ and $T_{X}(-\log \,\ddd)$ are  locally free, furthermore, the Poincar\'e residue map
$$
\rm{Res}: \Omega_X^1(\log\, \ddd) \longrightarrow \OO_{\ddd} \cong \bigoplus_{i=1}^N\OO_{\ddd_i}
$$
give the following exact sequence of sheaves on X
\begin{eqnarray}\label{seq1}
0 \longrightarrow \Omega_X^1 \longrightarrow \Omega_X^1(\log\, \ddd) \stackrel{Res}{\longrightarrow} \bigoplus_{i=1}^N\OO_{\ddd_i} \longrightarrow 0,
\end{eqnarray}
\noindent where $\Omega_X^1$ is the sheaf of holomorphics $1$-forms on $X$ and $\ddd_1,\ldots,\ddd_N$ are the irreducible components of $\ddd$.

Now, if  $\ddd$ is such that $\cod(\Sing(\ddd))>2$ then there exist the following exact sequence of sheaves on $X$ (see V. I. Dolgachev \cite{IVDnew}):
\begin{eqnarray}
0 \longrightarrow \Omega_X^1\longrightarrow \Omega_X^1(\log\, \ddd) \longrightarrow \OO_{\ddd} \longrightarrow 0,
\end{eqnarray}

On the projective space $\PP^n$, if $\ddd$ is a smooth   hypersurface,   then there exist the following exact sequence of sheaves  (see E. Angeline \cite{EA}):
\begin{eqnarray}\label{p2.07}
0 \longrightarrow T_{\PP^n}(-\log\, \ddd) \longrightarrow \OO_{\PP^n}(1)^{n+1} \longrightarrow \OO_{\PP^n}(k)\longrightarrow 0,
\end{eqnarray}
\noindent where $ k$ is the degree of $\ddd$.

\subsection{Singular one-dimensional holomorphic foliations }
\begin{defi} 
Let $X$ be a connected  complex manifold. An one-dimensional holomorphic foliation is
given by the following data:
\begin{itemize}
  \item[$i)$] an open covering $\mathcal{U}=\{U_{\alpha}\}$ of $X$;
  \item [$ii)$] for each $U_{\alpha}$ an holomorphic vector field $\zeta_\alpha$ ;
  \item [$iii)$]for every non-empty intersection, $U_{\alpha}\cap U_{\beta} \neq \emptyset $, a
        holomorphic function $$f_{\alpha\beta} \in \mathcal{O}_X^*(U_\alpha\cap U_\beta);$$
\end{itemize}
such that $\zeta_\alpha = f_{\alpha\beta}\zeta_\beta$ in $U_\alpha\cap U_\beta$ and $f_{\alpha\beta}f_{\beta\gamma} = f_{\alpha\gamma}$ in $U_\alpha\cap U_\beta\cap U_\gamma$.
\end{defi}
We denote by $K_{\F}$ the line bundle defined by the cocycle $\{f_{\alpha\beta}\}\in \mathrm{H}^1(X, \mathcal{O}^*)$. Thus, a one-dimensional holomorphic  foliation $\F$ on $X$ induces  a global holomorphic section $\zeta_{\F}\in \mathrm{H}^0(X,T_X\otimes K_{\F})$.

The line bundle $T_{\F}:= (K_{\F})^* \hookrightarrow T_X$ is called the  \emph{tangente  bundle} of $\F$. The singular set of $\F$ 
is $\singf=\{\zeta_{\F}=0\}$. We will assume that $\cod(\singf)\geq 2$.

\begin{defi} 
Let $V $ an analytic subspace of  a complex manifold $X$.  We say that $V$ is invariant by a foliation $\F$ if $T_{\F}|_{V}\subset (\Omega^1_V)^*$.
If $V$ is a hypersurface we say that $\F$ is \textit{logatithmic along $V$}.
\end{defi}

\begin{defi} 
A foliation on a  complex projective space $\mathbb{P}^n$  is  called by \textit{projective foliation}.
Let $\F $ be  a   projective foliation   with tangent bundle $T_{\F}=\mathcal{O}_{\mathbb{P}^n}(r)$. The integer number $d:=r+1$ is called by the \textit{degree} of $\F$.
\end{defi}

\subsection{The GSV-Index}

X. Gomez-Mont, J. Sead and A. Verjovsky \cite{GSV} introduced the {\it GSV-index} for a holomorphic vector field over an analytic hypersurface, with isolated singularities, on a complex manifold,  generalizing the (classical) Poincar\'e-Hopf  index. The concept of GSV-index was extended to holomorphic vector field on more general contexts. For example, J. Seade e T. Suwa in \cite{SeaSuw1}, defined the GSV-index for holomorphic vector field on analytic subvariety type isolated complete intersection singularity. J.-P. Brasselet, J. Seade and T. Suwa in \cite{a03}, extended the notion of GSV-index for vector fields defined in certain types of analytical subvariety with non-isolated singularities. 

\hyphenation{using}
In \cite{a08}   X. Gomez-Mont  defined the {\it homological index} of holomorphic vector field on an analytic hypersurface with isolated singularities, which coincides with GSV-index. There is also the {\it virtual index}, introduced by D. Lehmann, M. Soares and T. Suwa \cite{a07}, that via Chern-Weil theory can be interpreted as the GSV-index. 
M. Brunella \cite{a09} also present the GSV-index for foliations on complex surfaces  by a different approach.

\hyphenation{sin-gu-la-ri-ties}
\hyphenation{sin-gu-la-ri-ty}

Let $X$ be an $n$-dimensional complex manifold, $\ddd$ a isolated hypersurface singularity on $X$ and let $\fol$ be a foliation on $X$ of dimension one, with isolated singularities. Suppose $\fol$ logarithmic along $\ddd$, i.e., the analytic hypersurface $\ddd$ is invariant by each holomorphic vector field that is a local representative of $\fol$. The GSV-index of $\fol$ in $x\in \ddd$ will be denoted by $GSV(\fol,\ddd,x)$. For definition and details on the GSV-index we refer to \cite{BarSaeSuw} and \cite{Suw2}.

\hyphenation{fo-lia-tion}

\subsection{The Logarithmic Index} \label{seclog} Recently,  A. G. Aleksandrov introduced in \cite{Ale_p} the notion of logarithmic index for logarithmic   a vector field  . Let $\fol$ be an one-dimensional holomorphic foliation on $X$ with isolated singularities and logarithmic along $\ddd$. Fixed a point $x\in X$, let $v\in T_{X}(-\log \,\ddd)|_{U}$ a germ of vector field on $(U,x)$ tangent to $\fol$. The interior multiplication $i_{v}$ induces the complex of logarithmic differential forms
$$
0 \longrightarrow \Omega^n_{X,x}(\log\, \ddd) \stackrel{i_v}{\longrightarrow} \Omega^{n-1}_{X,x}(\log\, \ddd)\stackrel{i_v}{\longrightarrow} \cdots\stackrel{i_v}{\longrightarrow} \Omega^{1}_{X,x}(\log\, \ddd)\stackrel{i_v}{\longrightarrow} \OO_{n,x}.
$$
Since all singularities of $v$ are isolated, the $i_v$-homology groups of the complex $\Omega^{\bullet}_{X}(\log\, \ddd)$ are finite-dimensional vector spaces (see \cite{Ale_p}). Thus, the Euler characteristic
$$
\chi (\Omega^{\bullet}_{X}(\log\, \ddd), i_v) = \displaystyle\sum_{i=0}^n(-1)^i dim H_i(\Omega^\bullet_{X,x} (\log\, D),i_v).
$$
\noindent of the complex of logarithmic differential forms is well defined.  Since this number does not depend on local representative $v$ of the foliation $\fol$ at the point $x$, we define the {\it logarithmic} index of $\fol$ at the point $x$ by
$$
Ind_{log\,D,x}(\fol):= \chi (\Omega^{\bullet}_{X}(\log\, \ddd), i_v).
$$
\noindent It follows from the definition that $Ind_{log\,D,x}(\fol) = 0$ for all $x\in X - \singf$.  

We have the following important property (see \cite{Ale_p}):  

\begin{prop}  \cite{Ale_p} \label{prop004}
Let $X$, $\ddd$ and $\fol$ be as described above. Then, for each $x\in \Sing(\fol)\cap \ddd$ we have 
\begin{eqnarray}\label{eq001}
Ind_{log\,\ddd,x}(\fol) = \mu_x(\fol) - GSV(\fol,\ddd,x)
\end{eqnarray}
 where $\mu_x(v)$ and $GSV(v,\ddd,x)$ denote, respectively, the Milnor number and GSV index of $v$.
\end{prop}

\begin{remark}
If  $x\in \Sing(\fol)\cap \ddd_{reg}$, we obtain
\begin{eqnarray}\nonumber
Ind_{log\,\ddd,x}(\fol) = \mu_x(\fol) - \mu_x(\fol|_{{\ddd}_{reg}}),
\end{eqnarray}

\noindent  since, in this case, the GSV index of $\fol$ in $x$ coincides with the Milnor number of $\fol|_{{\ddd}_{reg}}$ in $x$. In particular, 
\begin{eqnarray}\nonumber
Ind_{log\,\ddd,x}(\fol)= 0,
\end{eqnarray}
\noindent whenever $x$ is a non-degenerate singularity of $\fol$. 
\end{remark}

\section{Proof of Theorem \ref{teo1_1}} \label{sec03}
To prove the Theorem \ref{teo1_1} we will firstly   prove the following result. 

\begin{teo}\label{teo_pre}
Let $X$ be an $n$-dimensional complex  compact manifold and  $\ddd$ a  smooth  hypersurface on $X$.  Then for all line bundle $L$  on $X$  we have 
\begin{eqnarray}\label{teo002}
\int_X c_{n}(T_X(-\log\, \ddd)- L) = \int_X c_{n}(T_X- L) - \int_{\ddd} c_{n-1}(T_X - [\ddd] - L).
\end{eqnarray}
\end{teo}
\begin{demostracion} By using   properties  of  Chern class we get
\begin{eqnarray}\nonumber
\displaystyle\int_{X}c_{n}(T_X(-\log\,  \ddd)- L) &=& \sum^{n}_{j=0}\int_{X}c_{n-j}(T_X(-\log\,  \ddd))c_1( L^*)^j\\\nonumber\\ \nonumber\\ \nonumber &=& \sum^{n}_{j=0}(-1)^{n-j}\int_{X}c_{n-j}(\Omega^1_X(\log\,  \ddd))c_1( L^*)^j.
\end{eqnarray}
\noindent On the one hand,  since   $\ddd$ is smooth     we can use  the exact sequence (\ref{seq1}) to  obtain 
\begin{eqnarray} \nonumber
c_i(\Omega_X^1(\log \ddd)) = \sum^{i}_{k=0}c_{i-k}(\Omega^1_X)c_k(\OO_{\ddd}), \forall i\in\{1,\ldots,n\}.
\end{eqnarray}
\noindent On the other hand, the Chern classes of $\OO_{\ddd}$ are
\begin{eqnarray} \label{r001}
c_k(\OO_{\ddd}) = c_k(\OO_X -  \OO (-\ddd) )= c_1([\ddd])^k, \,\,\,\, k =1,\ldots,n.  
\end{eqnarray}
 Thus,
$$
\int_{X}c_{n}(T_X(-\log\,  \ddd)- L) = \sum^{n}_{j=0}(-1)^{n-j}\int_{X}\left[\sum^{n-j}_{k=0}c_{n-j -k}(\Omega^1_X)c_1( [ \ddd])^k \right]c_1(L^*)^j.
$$
Now,  we  split  this sum in two parts  as follows:
\begin{eqnarray}\nonumber
\sum^{n}_{j=0}(-1)^{n-j}\int_{X}\left[\sum^{n-j}_{k=0}c_{n-j -k}(\Omega^1_X)c_1( [ \ddd])^k \right]c_1(L^*)^j = \\\nonumber \\\nonumber \displaystyle\sum^{n-1}_{j=0}(-1)^{n-j}\int_{X}\left[\sum^{n-j}_{k=1}c_{n-j -k}(\Omega^1_X)c_1( [ \ddd])^k \right ]c_1( L^*)^j +  \displaystyle\sum^{n}_{j=0}(-1)^{n-j}\int_{X}c_{n-j}(\Omega^1_X)c_1( L^*)^j.
\end{eqnarray}
\noindent In the first part  appears  all terms with $k\geq 1$ and in second one part  are the terms with $k=0$. By using the  Poincar\'e duality, we compute the first part as follow:
\begin{eqnarray}\nonumber
&&\displaystyle\sum^{n-1}_{j=0}(-1)^{n-j}\hspace*{-0.1cm}\int_{X}\hspace*{-0.1cm}\left[\sum^{n-j}_{k=1}c_{n-j -k}(\Omega^1_X)c_1([ \ddd])^k\hspace*{-0.08cm}\right]\hspace*{-0.08cm}c_1( L^*)^j =\\\nonumber\\\nonumber &=& \displaystyle\sum^{n-1}_{j=0}(-1)^{n-j}\hspace*{-0.08cm}\int_{ \ddd}\hspace*{-0.1cm}\left[\sum^{n-j}_{k=1}c_{n-j -k}(\Omega^1_X)c_1( [ \ddd])^{k-1}\hspace*{-0.08cm}\right]\hspace*{-0.08cm}c_1( L^*)^j=\\\nonumber\\\nonumber &=&
- \displaystyle\sum^{n-1}_{j=0}\int_{ \ddd}\left[\sum^{n-j}_{k=1}(-1)^{n-j-k}c_{n-j -k}(\Omega^1_X)(-1)^{k-1}c_1( [ \ddd])^{k-1}\right]c_1( L^{\ast})^j =
\\\nonumber\\\nonumber &=&
- \sum^{n-1}_{j=0} \displaystyle\int_{ \ddd} \left[\sum^{n-j}_{k=1}c_{n-j -k}(T_X)c_1( [ \ddd]^{\ast})^{k-1}\right]c_1(L^*)^j = \\\nonumber\\\nonumber &=&
- \sum^{n-1}_{j=0} \displaystyle\int_{ \ddd} c_{n-1-j}(T_X - [ \ddd])c_1(L^*)^j= \\\nonumber\\\nonumber &=&
- \displaystyle\int_{ \ddd} c_{n-1}(T_X - [ \ddd] - L).
\end{eqnarray}

Now, by basics proprieties of Chern classes   we compute the second sum as follow:
\begin{eqnarray}\nonumber
\displaystyle\sum^{n}_{j=0}(-1)^{n-j}\int_{X}c_{n-j}(\Omega^1_X)c_1( L^{\ast})^j &=&\displaystyle\sum^{n}_{j=0}(-1)^{n-j}\int_{X}(-1)^{n-j}c_{n-j}(T_X)c_1( L^{\ast})^j\\\nonumber
&=&\displaystyle\sum^{n}_{j=0}\int_{X}c_{n-j}(T_X)c_1( L^{\ast})^j\\\nonumber
&=& \displaystyle\int_{X}c_n(T_X - L).
\end{eqnarray}
\noindent Finally, we conclude that 
\begin{eqnarray}\nonumber
\sum^{n}_{j=0}(-1)^{n-j}\int_{X}\left[\sum^{n-j}_{k=0}c_{n-j -k}(\Omega^1_X)c_1( [ \ddd])^k \right]c_1(L^*)^j = - \displaystyle\int_{ \ddd} c_{n-1}(T_X - [ \ddd] - L) + \displaystyle\int_{X}c_n(T_X - L).
\end{eqnarray}
\noindent and this  proves the result. $\square$
\end{demostracion}

Now, we will prove the Theorem \ref{teo1_1}:

\begin{demostracion} Since   $\ddd$ is smooth,   we can invoke the formula (\ref{teo002}) of the  Theorem \ref{teo_pre}  to obtain the following equality

\begin{eqnarray}\nonumber
\displaystyle \int_{X} c_{n}(T_X(-\log\, \ddd)- T_{\fol}) = \int_{X} c_n(T_X - T_{\fol}) - \displaystyle\int_{\ddd} c_{n-1}(T_X - [\ddd]- T_{\fol}).
\end{eqnarray}

\noindent By hypothesis, the one-dimensional foliation $\fol$ is logarithmic along $\ddd$ and has only isolated singularities, then it follows from \cite{Suw2} that  the top Chern number of restriction $(T_X - [\ddd]- T_{\fol})|_{\ddd}$ coincides with the sum of GSV-Index of $\fol$ along $\ddd$. That is
$$
\displaystyle\int_{\ddd} c_{n-1}(T_X - [\ddd]- T_{\fol}) = \sum_{p\in  \Sing(\fol) \cap \ddd } GSV(\fol,\ddd,p),
$$
Hence
\begin{eqnarray}\nonumber
\displaystyle \int_{X} c_{n}(T_X(-\log\, \ddd)- T_{\fol}) &=& \int_{X} c_n(T_X - T_{\fol}) -\sum_{p\in  \Sing(\fol) \cap \ddd } GSV(\fol,\ddd,p)\\\nonumber \\\nonumber &=& \sum_{p\in \Sing(\fol)}\mu_p(\fol) - \sum_{p\in  \Sing(\fol) \cap \ddd } GSV(\fol,\ddd,p),
\end{eqnarray}
\noindent where in the last step we are using  the Baum-Bott classical formula (\ref{bbformula}).
Now, since $Ind_{log\,\ddd,p}(\fol)=0$ for all $p \in \Sing(\fol)\cap \ddd$, by Proposition \ref{prop004} we get  the following relation
$$
GSV(\fol,\ddd,p)=\mu_p(\fol),\   \forall\  p\in \Sing(\fol)\cap \ddd. 
$$
Therefore,  we obtain
\begin{eqnarray}\nonumber
\sum_{p\in \Sing(\fol)}\mu_p(\fol) - \sum_{p\in  \Sing(\fol) \cap \ddd } GSV(\fol,\ddd,p) = \sum_{p\in \Sing(\fol)\cap (X\setminus \ddd)}\mu_p(\fol) 
\end{eqnarray}

\noindent and the desired formula is proved. $\square$

\end{demostracion}

\section{Proof of Theorem 2} \label{sec05}
In this section we will consider $\ddd = \ddd_ {1} \cup \ldots \cup \ddd_{N}$ an analytic hypersurface on $X$,  with normal crossing singularities. Fixing an irreducible component, say $\ddd_N$, we define
$$
\hat{\ddd}_N:= \bigcup^{N-1}_{j=1} \ddd_{j} \,\,\,\,\,\,\,\,\,\, and \,\,\,\,\,\,\,\,\,\, 
\hat{\ddd}_N|\ddd_N:= \bigcup^{N-1}_{j=1} \ddd_j\cap\ddd_N.
$$
We note that $\hat{\ddd}_N|\ddd_N$ is an analytic hypersurface on $\ddd_N$ with normal crossings singularities and $N-1$ irreducible components. We will use the following multiple index notation: for each multi-index $J = (j_1,\ldots,j_N)$ and $J'=(j_1',\ldots,j_{N-1}')$, with $1\leq j_l,j'_t\leq n$, we denote 
\begin{eqnarray}\nonumber
c_1(\ddd)^J & = & c_1([\ddd_1])^{j_1}\cdots c_1([\ddd_N])^{j_N}, \\\nonumber
c_1(\hat{\ddd}_N)^{J'} &=& c_1([\ddd_1])^{j'_1}\cdots c_1([\ddd_{N-1}])^{j'_{N-1}}.
\end{eqnarray}

\begin{lema} \label{p13} In the above conditions, for each $i = 1,\ldots, n$, we have
\begin{eqnarray}\nonumber
\displaystyle c_{i}(\Omega_{X}^1(\log\, \ddd)) = \sum^{i}_{k=0}\sum_{\mid J \mid = k} c_{i-k}(\Omega_X^1)c_1(\ddd)^{J}.
\end{eqnarray}
\end{lema}
\begin{demostracion} Since $\ddd$ is an analytic hypersurface with normal crossing singularities, the Poincar\'e residue map
$$
Res: \Omega_X^1(\log\, \ddd) \longrightarrow \OO_{\ddd} \cong \bigoplus_{i=1}^N\OO_{\ddd_i}
$$
induces  the following exact sequence 
\begin{eqnarray}\nonumber
0 \longrightarrow \Omega_X^1 \longrightarrow \Omega_X^1(\log\, \ddd) \stackrel{Res}{\longrightarrow} \bigoplus_{i=1}^N\OO_{\ddd_i} \longrightarrow 0.
\end{eqnarray}
\noindent By  using this exact sequence, we get 
\begin{eqnarray} \nonumber
c_i(\Omega_X^1(\log\, \ddd)) &=& \sum^{i}_{k=0}c_{i-k}(\Omega^1_X)c_k(\bigoplus_{i=1}^N\OO_{\ddd_i}) \\\nonumber \\\nonumber  &=& \sum^{i}_{k=0}c_{i-k}(\Omega^1_X)\left (\sum_{j_1+\ldots +j_N=k}c_{j_1}(\OO_{\ddd_1})\ldots c_{j_N}(\OO_{\ddd_N})\right) \\\nonumber \\\nonumber &=&\sum^{i}_{k=0}c_{i-k}(\Omega_X^1)\left(\sum_{j_1+\ldots+j_N = k}c_1([\ddd_1])^{j_1}\ldots c_1([\ddd_N])^{j_N}\right),
\end{eqnarray}
\noindent where in last equality we use  the following relations
$$
c_i(\OO_{\ddd_j}) = c_1([\ddd_j])^i,\,\,\,\, i = 1,\ldots,n,
$$
 which can   be obtained of (\ref{r001}). $\square$
\end{demostracion}

\begin{lema} \label{p1300} In the above conditions, for each $i = 1,\ldots, n-1$,
\begin{eqnarray}\nonumber
\displaystyle c_{i}(\Omega_{X}^1)|_{\ddd_N} = c_{i}(\Omega_{\ddd_N}^1) - c_{i-1}(\Omega_{\ddd_N}^1)c_i([\ddd_N])|_{\ddd_N}.
\end{eqnarray}
\end{lema}
\begin{demostracion} 
It follows from by taking the total  Chern class in the exact sequence
$$
0\to  T_{\ddd_N}  \to T_X|_{\ddd_N}  \to  [\ddd_N]|_{\ddd_N}\to 0.
$$
\end{demostracion}

\begin{lema} \label{p0013} In the above conditions, if $L$ is a holomorphic line bundle  on $X$, then
the following relations hold: 
\begin{eqnarray} \label{rel001}
\int_{X}c_n(T_X(-\log\, \ddd) -  L) =\\\nonumber \sum_{j=0}^n\sum^{n-j}_{k=0} \sum_{\mid J \mid = k}  \int_{X} (-1)^{n-j}c_{n-j-k}(\Omega_X^1) c_1(\ddd)^{J}c_1(L^{\ast})^j.
\end{eqnarray}
\noindent In particular,
\begin{eqnarray} \label{rel002}
\int_{X}c_n(T_X(-\log\, \hat{\ddd}_N) -  L) =\\\nonumber \sum_{j=0}^n\sum^{n-j}_{k=0} \sum_{\mid J' \mid = k}  \int_{X} (-1)^{n-j}c_{n-j-k}(\Omega_X^1) c_1(\hat{\ddd}_N)^{J'}c_1(L^{\ast})^j.
\end{eqnarray}
\noindent and
\begin{eqnarray} \label{rel003}
\int_{\ddd_N}c_{n-1}(T_{\ddd_N}(-\log\, (\hat{\ddd}_N|\ddd_N)) - L|_{\ddd_N}) = \\\nonumber \sum_{j=0}^{n-1}\sum^{n-1-j}_{k=0} \sum_{\mid J' \mid = k}  \int_{\ddd_N} (-1)^{n-1-j}c_{n-1-j-k}(\Omega_{\ddd_N}^1) c_1(\hat{\ddd}_N)^{J'}c_1(L^{\ast})^j.
\end{eqnarray}
\end{lema}

\begin{demostracion} By  basics proprieties of Chern classes, we get
\begin{eqnarray}\nonumber
\displaystyle \int_{X}c_n(T_X(-\log\, \ddd) - L) &=& \int_{X} \sum_{j=0}^n c_{n-j}(T_{X}(-\log\, \ddd))c_1(L^{\ast})^j\\\nonumber
&=& \int_{X}\sum_{j=0}^n(-1)^{n-j} c_{n-j}(\Omega_{X}^1(\log\, \ddd))c_1(L^{\ast})^j.
\end{eqnarray}
\noindent By Lemma \ref{p13}, we get
$$
c_{n-j}(\Omega_{X}^1(\log \, \ddd)) = \sum^{n-j}_{k=0} \sum_{\mid J \mid = k} c_{n-j-k}(\Omega_X^1)c_1(\ddd)^{J}.
$$ 
 Substituting this, we obtain (\ref{rel001}).
The relation (\ref{rel002}) is obtained by taking $\ddd=\hat{\ddd}_N$ in relation (\ref{rel001}). Analogously, applying  the relation (\ref{rel002}), we can obtain (\ref{rel003}) by taking $X = \ddd_N$ as a complex manifold of dimension $n-1$ and $\ddd = \hat{D}_N | D_N$ as an analytic subvariety of  $\ddd_N$ with normal crossings. $\square$
\end{demostracion}

\begin{prop}\label{lema4}
In the above conditions, if $L$  is a holomorphic line bundle  on $X$, then
{\small {\begin{eqnarray}\nonumber
\int_{X}c_n(T_X(-\log\, \ddd) -  L) &=& \int_{X}c_n(T_X(-\log (\hat{\ddd}_N) - L) - \int_{\ddd_N}c_{n-1}(T_{\ddd_N} 
(- \log( \hat{\ddd}_N|\ddd_N) - L|_{\ddd_N}).
\end{eqnarray}}}
\end{prop}

\begin{demostracion} By Lemma \ref{p0013}, it is sufficient to show that the following equality occurs

\begin{eqnarray}\nonumber
\sum_{j=0}^n\sum^{n-j}_{k=0} \sum_{\mid J \mid = k}  \int_{X} (-1)^{n-j}c_{n-j-k}(\Omega_X^1) c_1(\ddd)^{J}c_1(L^{\ast})^j &=&\\\nonumber 
\sum_{j=0}^n\sum^{n-j}_{k=0} \sum_{\mid J' \mid = k}  \int_{X} (-1)^{n-j}c_{n-j-k}(\Omega_X^1) c_1(\hat{\ddd}_N)^{J'}c_1(L^{\ast})^j &-&\\\nonumber - \sum_{j=0}^{n-1}\sum^{n-1-j}_{k=0} \sum_{\mid J' \mid = k}  \int_{\ddd_N} (-1)^{n-1-j}c_{n-1-j-k}(\Omega_{\ddd_N}^1) c_1(\hat{\ddd}_N)^{J'}c_1(L^{\ast})^j.
\end{eqnarray}

\noindent Indeed, we can decompose the sum on the left hand side into the terms with $k=0$ and those with $k\geq 1$ as follows:

\begin{eqnarray}\label{eq00003}
\sum_{j=0}^n\sum^{n-j}_{k=0} \sum_{\mid J \mid = k}  \int_{X} (-1)^{n-j}c_{n-j-k}(\Omega_X^1) c_1(\ddd)^{J}c_1(L^{\ast})^j &=&  \\\nonumber\\\nonumber\\\nonumber  =  \displaystyle \sum^{n}_{j=0}\int_{X}(-1)^{n-j}c_{n-j}(\Omega^1_X)c_1(L^{\ast})^j + \sum_{j=0}^{n-1}\sum^{n-j}_{k=1} \sum_{\mid J \mid = k}  \int_{X} (-1)^{n-j}c_{n-j-k}(\Omega_X^1) c_1(\ddd)^{J}c_1(L^{\ast})^j.
\end{eqnarray}
\noindent The second sum on the right hand side can readily be computed. In fact, 
\begin{eqnarray}\nonumber
\sum_{j=0}^{n-1}\sum^{n-j}_{k=1} \sum_{\mid J \mid = k}  \int_{X} (-1)^{n-j}c_{n-j-k}(\Omega_X^1) c_1(\ddd)^{J}c_1(L^{\ast})^j = \\\nonumber\\\nonumber
\sum_{j=0}^{n-1}\sum^{n-j}_{k=1} \sum_{\mid J' \mid = k}  \int_{X} (-1)^{n-j}c_{n-j-k}(\Omega_X^1) c_1(\hat{\ddd}_N)^{J'}c_1(L^{\ast})^j + \\\nonumber\\\nonumber  + \sum_{j=0}^{n-1}\sum^{n-j}_{k=1} \sum_{\substack{\mid J \mid = k \\ j_N \geq 1}}  \int_{X} (-1)^{n-j} c_{n-j-k}(\Omega_X^1) c_1([\ddd_1])^{j_1}\ldots c_1([\ddd_N])^{j_N}c_1(L^{\ast})^j.
\end{eqnarray}
\noindent  By using the fact that $c_1 ([\ddd_N])$ is Poincar\'e dual to the fundamental class of $\ddd_N$,
we obtain:
\begin{eqnarray}\nonumber
\sum_{j=0}^{n-1}\sum^{n-j}_{k=1} \sum_{\mid J \mid = k}  \int_{X} (-1)^{n-j}c_{n-j-k}(\Omega_X^1) c_1(\ddd)^{J}c_1(L^{\ast})^j = \\\nonumber\\\nonumber
\sum_{j=0}^{n-1}\sum^{n-j}_{k=1} \sum_{\mid J' \mid = k}  \int_{X} (-1)^{n-j}c_{n-j-k}(\Omega_X^1) c_1(\hat{\ddd}_N)^{J'}c_1(L^{\ast})^j + \\\nonumber\\\nonumber + \sum_{j=0}^{n-1}\sum^{n-j}_{k=1} \sum_{\substack{\mid J \mid = k \\ j_N \geq 1}}  \int_{\ddd_N} (-1)^{n-j} c_{n-j-k}(\Omega_X^1) c_1([\ddd_1])^{j_1}\ldots c_1([\ddd_N])^{j_N - 1}c_1(L^{\ast})^j.
\end{eqnarray}

\noindent Now, using the relation of Lemma \ref{p1300}, we get

\begin{eqnarray}\nonumber
\sum_{j=0}^{n-1}\sum^{n-j}_{k=1} \sum_{\substack{\mid J \mid = k \\ j_N \geq 1}}  \int_{\ddd_N} (-1)^{n-j} c_{n-j-k}(\Omega_X^1) c_1([\ddd_1])^{j_1}\ldots c_1([\ddd_N])^{j_N - 1}c_1(L^{\ast})^j = \\\nonumber\\\nonumber = \sum_{j=0}^{n-1}\sum^{n-j}_{k=1} \sum_{\substack{\mid J \mid = k \\ j_N \geq 1}}  \int_{\ddd_N} (-1)^{n-j} c_{n-j-k}(\Omega_{\ddd_N}^1) c_1([\ddd_1])^{j_1}\ldots c_1([\ddd_N])^{j_N - 1}c_1(L^{\ast})^j - \\\nonumber\\\nonumber  - \sum_{j=0}^{n-1}\sum^{n-j-1}_{k=1} \sum_{\substack{\mid J \mid = k \\ j_N \geq 1}}  \int_{\ddd_N} (-1)^{n-j} c_{n-j-1-k}(\Omega_{\ddd_N}^1) 
c_1(\ddd)^{J}c_1(L^{\ast})^j = \\\nonumber\\\nonumber  = \sum_{j=0}^{n-1}  \int_{\ddd_N} (-1)^{n-j} c_{n-j-1}(\Omega_{\ddd_N}^1) c_1(L^{\ast})^j + \sum_{j=0}^{n-1}\sum^{n-j-1}_{k=1} \sum_{\mid J' \mid = k}  \int_{\ddd_N} (-1)^{n-j} c_{n-j-1-k}(\Omega_{\ddd_N}^1) c_1(\hat {\ddd}_N)^{J'}c_1(L^{\ast})^j = \\\nonumber\\\nonumber = - \sum_{j=0}^{n-1}\sum^{n-1-j}_{k=0} \sum_{\mid J' \mid = k}  \int_{\ddd_N} (-1)^{n-1-j} c_{n-1-j-k}(\Omega_{\ddd_N}^1) c_1(\hat {\ddd}_N)^{J'}c_1(L^{\ast})^j. 
\end{eqnarray}

Hence,

\begin{eqnarray}\nonumber
\sum_{j=0}^{n-1}\sum^{n-j}_{k=1} \sum_{\mid J \mid = k}  \int_{X} (-1)^{n-j}c_{n-j-k}(\Omega_X^1) c_1(\ddd)^{J}c_1(L^{\ast})^j = \\\nonumber\\\nonumber
\sum_{j=0}^{n-1}\sum^{n-j}_{k=1} \sum_{\mid J' \mid = k}  \int_{X} (-1)^{n-j}c_{n-j-k}(\Omega_X^1) c_1(\hat{\ddd}_N)^{J'}c_1(L^{\ast})^j - \\\nonumber\\\nonumber - \sum_{j=0}^{n-1}\sum^{n-1-j}_{k=0} \sum_{\mid J' \mid = k}  \int_{\ddd_N} (-1)^{n-1-j} c_{n-1-j-k}(\Omega_{\ddd_N}^1) c_1(\hat {\ddd}_N)^{J'}c_1(L^{\ast})^j,
\end{eqnarray}

\noindent and we completed the calculation of the second sum.

Replacing it in the initial equality (\ref{eq00003}), we obtain

\begin{eqnarray}\nonumber 
\sum_{j=0}^n\sum^{n-j}_{k=0} \sum_{\mid J \mid = k}  \int_{X} (-1)^{n-j}c_{n-j-k}(\Omega_X^1) c_1(\ddd)^{J}c_1(L^{\ast})^j &=&  \\\nonumber\\\nonumber  =  \displaystyle \sum^{n}_{j=0}\int_{X}(-1)^{n-j}c_{n-j}(\Omega^1_X)c_1(L^{\ast})^j + \sum_{j=0}^{n-1}\sum^{n-j}_{k=1} \sum_{\mid J' \mid = k}  \int_{X} (-1)^{n-j}c_{n-j-k}(\Omega_X^1) c_1(\hat{\ddd}_N)^{J'}c_1(L^{\ast})^j - \\\nonumber\\\nonumber - \sum_{j=0}^{n-1}\sum^{n-1-j}_{k=0} \sum_{\mid J' \mid = k}  \int_{\ddd_N} (-1)^{n-1-j} c_{n-1-j-k}(\Omega_{\ddd_N}^1) c_1(\hat{\ddd}_N)^{J'}c_1(L^{\ast})^j = \\\nonumber\\\nonumber  = \sum_{j=0}^n\sum^{n-j}_{k=0} \sum_{\mid J' \mid = k}  \int_{X} (-1)^{n-j}c_{n-j-k}(\Omega_X^1) c_1(\hat{\ddd}_N)^{J'}c_1(L^{\ast})^j - \\\nonumber\\\nonumber - \sum_{j=0}^{n-1}\sum^{n-1-j}_{k=0} \sum_{\mid J' \mid = k}  \int_{\ddd_N} (-1)^{n-1-j}c_{n-1-j-k}(\Omega_{\ddd_N}^1) c_1(\hat{\ddd}_N)^{J'}c_1(L^{\ast})^j .
\end{eqnarray}

$\square$

\end{demostracion}

Now, we will prove the Theorem $2$:
\\
\begin{demostracion} We will prove by induction on the number of irreducible components of $\ddd$. Indeed, if the number of irreducible component of $\ddd$ is 1, then $\ddd$ is smooth. By hypothesis, the singularities of $\fol$ are non-degenerate, and thus the theorem follows  from Theorem \ref{teo1_1}.

\hyphenation{hy-po-the-sis}

Let us suppose that for every analytic hypersurface on $X$, satisfying the hypothesis of theorem and having $N-1$ irreducible components, the formula (\ref{pp1}) holds. Let $\ddd$ be an analytic hypersurface on $X$ with $N$ irreducible components, satisfying the hypotheses of the theorem. We will prove that the formula (\ref{pp1}) is true for $\ddd$.

We know that $\hat{\ddd}_N$ is an analytic hypersurface on $X$ and  $\hat{\ddd}_N|{\ddd}_N$ is an analytic hypersurface on $\ddd_N$, both with normal crossing singularities and having exactly $N - 1$ irreducible components. Moreover, $\fol$ and its restriction  $\fol|_{\ddd_N}$ on $\ddd_N$ are logarithmic along $\ddd_N$ and $\hat{\ddd}_N|{\ddd}_N$, respectively. Thus, we can use the induction hypothesis and we obtain
\begin{eqnarray}\label{expr108}
 \sum_{p\in Sing \left (\fol \right )\cap \left (X\setminus \hat{\ddd}_N \right)}\mu_p(\fol) = \int_{X}c_n(T_X(-\log\, \hat{\ddd}_N)- T_{\fol})
\end{eqnarray}

\noindent and

\begin{eqnarray}\label{expr109}
\sum_{p\in Sing \left(\fol\right) \cap \left [\ddd_{N}\setminus \left(\hat{\ddd}_N|\ddd_{N}\right)\right]}\mu_p(\fol)  =  \int_{\ddd_{N}}c_{n-1}(T_{\ddd_N}(-\log\,(\hat{\ddd}_N|\ddd_{N}))- T_{\fol}|_{\ddd_{N}}).
\end{eqnarray}
\noindent By using the following identity
$
X -  \ddd = (X - \hat{\ddd}_N) - [\ddd_N - (\hat{\ddd}_N \cap \ddd_N)], 
$
we get 
\begin{eqnarray}\nonumber
\sum_{p\in Sing \left(\fol\right) \cap \left(X\setminus \ddd\right)}\mu_p(\fol) =  \sum_{p\in Sing \left(\fol \right) \cap \left (X\setminus \hat{\ddd}_N \right)} \mu_p(\fol) - \sum_{p\in Sing \left(\fol \right) \cap \left[\ddd_{N}\setminus \left(\hat{\ddd}_N|\ddd_N \right) \right]}\mu_p(\fol).
\end{eqnarray}
\noindent Therefore, by (\ref{expr108}) and (\ref{expr109}), we get
\begin{eqnarray}\nonumber
\sum_{p\in Sing \left(\fol\right) \cap \left(X\setminus \ddd\right)}\mu_p(\fol) =  \int_{X}c_n(T_X(-\log\, \hat{\ddd}_N)- T_{\fol}) - \int_{\ddd_{N}}c_{n-1}(T_{\ddd_N}(-\log\,(\hat{\ddd}_N|\ddd_{N}))- T_{\fol}|_{\ddd_{N}}),
\end{eqnarray}
\noindent  and   we obtain the desired equality by  applying   the Proposition \ref{lema4}. Thus, we prove that the formula (\ref{pp1}) is true for $\ddd$ and the proof of the theorem follows by induction. 
$\square$
\end{demostracion}

\section{Application: a Poincar\'e-Hopf  type Theorem} \label{sec04}

In this section  we will prove a  Poincar\'e-Hopf  type Theorem for non-compact complex manifolds.  More precisely, we prove the following:  

\medskip

\noindent {\bf Corollary \ref{cor1_1}.} {\it Let $\tilde{X}$ be an $n$-dimensional complex manifold such that $\tilde{X} = X - \ddd$, where $X$ is an $n$-dimensional complex compact  manifold, $\ddd$ is a reduced  normal crossing hypersurface on $X$.  Let $\fol$ be a foliation on $X$ of dimension one given by a global holomorphic vector field, with isolated singularities (non-degenerates) and logarithmic along $\ddd$. Then 
$$
 \chi(\tilde{X}) =\sum_{p\in \Sing(\fol)\cap  \tilde{X} }\mu_p(\fol),
$$
where $\mu_p(\fol)$ denotes the Milnor number of $\fol$ on $p$.}
\medskip

\begin{demostracion} On the one hand,  it follows from    Norimatsu-Silvotti -Aluffi  Theorem  that 
\begin{eqnarray}\nonumber
\displaystyle \int_{X} c_{n}(T_X(-\log\, \ddd)) =  \chi (\tilde{X}).
\end{eqnarray}
On the other hand, since   $\ddd$ is a  normal crossing   hypersurface, it follows from  Theorem $2$  that  
\begin{eqnarray}\nonumber
\displaystyle \int_{X} c_{n}(T_X(-\log\, \ddd)) &=& \sum_{p\in \Sing(\fol)\cap  \tilde{X} }\mu_p(\fol).
\end{eqnarray}
\noindent  This shows the result.
\end{demostracion}

\section{Application to one-dimensional projective foliations} \label{sec06}

In this section we give an optimal  description  for a smooth hypersurface $\ddd$ invariant by an one-dimensional foliation $\fol$ on  $\mathbb{P}^n$ satisfying $\Sing(\fol) \subsetneq \ddd.$ More precisely,   we will prove the Theorem $3$.
  
Firstly, we need some preliminary results.

\begin{lema}\label{lema001}
Let $f(x,y) = \displaystyle\sum^{n}_{i=0} \sum^{n-i}_{j=0}\binom{n + 1}{n-i-j}
x^jy^i$, with $n\in \NN$. 

\bigskip

\begin{enumerate}
\item[(i)] If $x\neq y$, then $f(x,y) = \displaystyle\frac{(1+x)^{n+1} - (1+y)^{n+1}}{x-y}$.
\bigskip

\item[(ii)] If $x = y$, then $f(x,y) = (n+1)(1+x)^n$
\end{enumerate}

\end{lema} \begin{demostracion}  Developing the summation in the following  triangular format 
\begin{equation}
\begin{array}{lllllllllclll}\nonumber
\displaystyle f(x,y) =  \displaystyle\binom{n+1}{n}   \hspace*{- 0.2 cm}&+&\hspace*{- 0.2 cm} \displaystyle\binom{n+1}{n-1}x \hspace*{- 0.2 cm}&+&\hspace*{- 0.2 cm} \displaystyle\binom{n+1}{n-2}x^2 \hspace*{- 0.2 cm}&+&\hspace*{- 0.2 cm} \hspace*{- 0.2 cm}&\ldots&\hspace*{- 0.2 cm} \hspace*{- 0.2 cm} &+& \hspace*{- 0.3 cm}\displaystyle\binom{n+1}{0}x^{n}\hspace*{- 0.2 cm}&+&\\
\nonumber\\ \nonumber

\hspace*{- 0.2 cm} &+& \hspace*{- 0.2 cm} \displaystyle\binom{n+1}{n-1}y \hspace*{- 0.2 cm}&+&\hspace*{- 0.2 cm} \displaystyle\binom{n+1}{n-2}yx \hspace*{- 0.2 cm}&+&\hspace*{- 0.2 cm}  \hspace*{- 0.2 cm} &\ldots&\hspace*{- 0.2 cm}  \hspace*{- 0.2 cm}&+&\hspace*{- 0.2 cm} \displaystyle\binom{n+1}{0}yx^{n-1} \hspace*{- 0.2 cm}&+&\hspace*{- 0.2 cm}\\
\nonumber\\ \nonumber

\hspace*{- 0.2 cm} && \hspace*{- 0.2 cm} \hspace*{- 0.2 cm}&+&\hspace*{- 0.2 cm} \displaystyle\binom{n+1}{n-2}y^2  \hspace*{- 0.2 cm}&+&\hspace*{- 0.2 cm}
\hspace*{- 0.2 cm} &\ldots&\hspace*{- 0.2 cm}    \hspace*{- 0.2 cm}&+&\hspace*{- 0.2 cm} \displaystyle\binom{n+1}{0}y^2x^{n-2} \hspace*{- 0.2 cm}&+& \hspace*{- 0.2 cm}
\\ \nonumber

\hspace*{- 0.2 cm} && \hspace*{- 0.2 cm} \hspace*{- 0.2 cm}&&\hspace*{- 0.2 cm} \hspace*{- 0.2 cm}&& \hspace*{- 0.2 cm}  \hspace*{- 0.2 cm}&& \hspace*{- 0.2 cm}  \hspace*{- 0.2 cm}&\vdots&\hspace*{- 0.2 cm}     \hspace*{- 0.2 cm}&&\hspace*{- 0.2 cm}\\
\nonumber

\hspace*{- 0.2 cm} && \hspace*{- 0.2 cm} \hspace*{- 0.2 cm}&&\hspace*{- 0.2 cm} \hspace*{- 0.2 cm}&&\hspace*{- 0.2 cm} \hspace*{- 0.2 cm}&&\hspace*{- 0.2 cm}  \hspace*{- 0.2 cm}&+&\hspace*{- 0.2 cm} \displaystyle\binom{n+1}{0}y^{n}.\hspace*{- 0.2 cm}&&
\end{array}
\end{equation}
\noindent We can put in evidence the common factor in   each columns and we get the following forma:
\begin{eqnarray}\label{eq001}
f(x,y) = \displaystyle \sum^n_{k=0}\left [\displaystyle\binom{n+1}{n-k}\left (\sum_{j=0}^k x^{k-j}y^{j}\right)\right].
\end{eqnarray}

\noindent (i)  Suppose $x\neq y$. Since
$$
x^k +x^{k-1}y +\ldots +xy^{k-1} +y^k = \displaystyle\frac{x^{k+1} - y^{k+1}}{x-y},\,\, 0\leq k\leq n,
$$
 we obtain
\begin{eqnarray}\nonumber
f(x,y) &=& \displaystyle \sum_{k = 0}^n \left [ \displaystyle \binom{n+1}{n-k}\left(\frac{x^{k+1} - y^{k+1}}{x-y}\right) \right]\\\nonumber\\\nonumber
&=& \frac{1}{x-y}\left[\sum_{k = 0}^n \binom{n+1}{n-k} x^{k+1} - \sum_{k = 0}^n \binom{n+1}{n-k} y^{k+1} \displaystyle \right]\\\nonumber\\\nonumber
&=& \displaystyle\frac{(1+x)^{n+1} - (1+y)^{n+1}}{x-y},
\end{eqnarray}
\noindent where in the last equality we have used the binomial theorem.\\

\noindent (ii) Consider the case where $x=y$. By equality (\ref{eq001}) we have 

\begin{eqnarray}\nonumber
f(x,y) &=& \binom{n+1}{n} + \binom{n+1}{n-1}2x + \binom{n+1}{n-1}3x^2 +\ldots+ \binom{n+1}{0}(n+1)x^n.
\end{eqnarray}

\noindent Hence, by using the binomial theorem we obtain

\begin{eqnarray}\nonumber
f(x,y) &=& \frac{d}{dx}\left[ (1+x)^{n+1}\right]\\\nonumber\\\nonumber
&=& (n+1)(1+x)^{n}.
\end{eqnarray}

$\square$\end{demostracion}

\begin{lema}\label{lemmm}
Let $k, d$ and $n$ natural numbers, with $k\geq 1$, $d \geq 0$ and $n\geq 2$. Consider the natural number $\delta(k,d,n)$ defined by the relation
$$
\delta(k,d,n) = \displaystyle\sum^{n}_{i=0} \sum^{n-i}_{j=0}\binom{n + 1}{n-i-j}
(-k)^j (d-1)^i.
$$

\vspace*{-0.13 cm}

\noindent Then $\delta(k,d,n)$ satisfies the following conditions:\\

\noindent {\bf (i)} If $n$ is odd, then:\\

{\bf(a)} $\delta(k,d,n)>0 \Longleftrightarrow k<d+1$;\\

{\bf(b)} $\delta(k,d,n)=0 \Longleftrightarrow k=d+1$;\\

{\bf(c)} $\delta(k,d,n)<0 \Longleftrightarrow k>d+1$.\\\\

\noindent {\bf(ii)} If $n$ is even, then $\delta(k,d,n) \geq 0$ and moreover:\\

{\bf(a)} $\delta(k,d,n)>0 \Longleftrightarrow \left \{ \begin{array}{lll}
  k\neq d+1\\
\mbox{or}\\
k=d+1,\,\,\mbox{with}\,\, d\neq 0 \\
\end{array}
\right.\,\,\,\,\,\,;$\\

{\bf(b)} $\delta(k,d,n) = 0 \Longleftrightarrow k=1$ and $d=0$.\\\\
\noindent {\bf(iii)} $\delta(k,d,n) = \displaystyle\sum^{n}_{i=0}(-1)^i(k-1)^id^{n-i}$.

\end{lema}

\begin{demostracion} In order to prove the present lemma, we can consider $x=-k$ and $y=d-1$ in $f(x,y)$ of the lemma \ref{lema001}. Hence, we obtain

\begin{eqnarray} \nonumber
\delta(k,d,n) = \displaystyle \frac{(1-k)^{n+1} - d^{n+1}}{-k-d+1},\,\,\,\mbox{if $k\neq 1$ or $d\neq 0$} 
\end{eqnarray}

\noindent and

\begin{eqnarray}\label{cond02}
\delta(k,d,n) = 0,\,\,\,\mbox{if $k=1$ and $d= 0$}.
\end{eqnarray}

\noindent The proof of itens (i) and (ii) can readily be obtained by the study of sign of the expression  
$$
\displaystyle\frac{(1-k)^{n+1} - d^{n+1}}{-k-d+1}
$$
 and also using the relation (\ref{cond02}).

Now, let us consider the  summation  $\displaystyle\sum^{n}_{i=0}(-1)^i(k-1)^id^{n-i}$. We have:
\begin{eqnarray}\nonumber
\sum^{n}_{i=0}(-1)^i(k-1)^id^{n-i}&=& d^n\left[\sum^{n}_{i=0}(-1)^i(k-1)^id^{-i}\right]\\\nonumber\\\nonumber
&=& d^n\left[\sum^{n}_{i=0}\left(\frac{(-1)(k-1)}{d}\right)^i\right].
\end{eqnarray}
\noindent By the  property
$$
\forall a\in\ZZ,\,\,\,\,\, 1 +a + a^2+\ldots +a^{n} = \displaystyle\frac{1 - a^{n+1}}{1-a}
$$
\noindent we get
\begin{eqnarray}\nonumber
\sum^{n}_{i=0}(-1)^i(k-1)^id^{n-i}&=& \displaystyle\frac{(k-1)^{n+1} - d^{n+1} }{-k-d+1}.
\end{eqnarray}

\noindent Hence, this proves  the equality of item (iii).

$\square$\end{demostracion}

\begin{lema}\label{lema0001}
Let $\ddd\subset \mathbb{P}^n$ a smooth and irreducible hypersurface of degree $  k$. Then, for $l = 1,\ldots,n$, we obtain 
\begin{eqnarray} \label{formula0001}
\displaystyle c_{l}(T_{\CIPT}(-\log \,\ddd)) = \displaystyle\left [\displaystyle\sum^{l}_{j=0} \binom{n + 1}{l-j }(-1)^jk^j\right ]c_1(\OO_{\CIPT}(1))^l.
\end{eqnarray}

\end{lema}

\begin{demostracion} The   formula (\ref{formula0001}) can   be obtained by considering the recursion
\begin{eqnarray}\nonumber
\displaystyle c_{j+1}(T_{\CIPT}(-\log\, \ddd)) = \displaystyle \binom{n + 1}{j+1}c_1(\OO_{\CIPT}(1))^{j+1} - c_{j}(T_{\CIPT}(-\log \,\ddd))(k\,c_1(\OO_{\CIPT}(1))),
\end{eqnarray}
\noindent  $j = 0,\ldots,n-1$, which can be obtained considering  the exact sequence (\ref{p2.07}):

\begin{eqnarray}\nonumber
0 \longrightarrow T_{\CIPT}(-\log\, \ddd) \longrightarrow \OO_{\CIPT}(1)^{n+1} \longrightarrow \OO_{\CIPT}(k)\longrightarrow 0.
\end{eqnarray}
$\square$
 \end{demostracion}

Now, we will prove the Theorem $3$:

\begin{demostracion} Let $deg(\ddd) = k$ and $deg(\fol) = d$. On the one hand, we have 
\begin{eqnarray}\nonumber
\displaystyle\int_{\CIPT}c_{n}(T_{\CIPT}(-\log \, \ddd)-T_{\fol}) = \displaystyle\sum^{n}_{i=0}\int_{\CIPT} c_{n-i}(T_{\CIPT}(-\log \, \ddd))c_1( T_{\fol}^*)^i.
\end{eqnarray}

\noindent Now, by using  the formula (\ref{formula0001}) in each $c_{n-i}(T_{\CIPT}(-\log \, \ddd))$, in the summation  above, we obtain
\begin{eqnarray}\nonumber
\displaystyle\int_{\CIPT}c_{n}(T_{\CIPT}(-\log \, \ddd)-T_{\fol}) = \displaystyle\sum^{n}_{i=0} \left [\sum^{n-i}_{j=0}\binom{n+1}{n-i-j} (-1)^jk^j\right ]\int_{\CIPT}c_1(\OO_{\CIPT}(1))^{n-i}c_1( T_{\fol}^*)^i.
\end{eqnarray}

\noindent On the other hand, the tangent bundle $T_{\fol}$ of foliation on $\CIPT$ is such that $T_{\fol}  = \OO_{\CIPT}(1-d)$. Therefore, we obtain $c_1(T_{\fol}^*) = (d-1)c_1(\OO_{\CIPT}(1))$. Hence,

\begin{eqnarray}\nonumber
\displaystyle\int_{\CIPT}c_{n}(T_{\CIPT}(-\log \, \ddd)-T_{\fol}) &=& \displaystyle\sum^{n}_{i=0} \left [\sum^{n-i}_{j=0} \binom{n + 1}{n-i-j}(-1)^jk^j\right ] (d-1)^i\int_{\CIPT}c_1(\OO_{\CIPT}(1))^{n}\\\nonumber\\\nonumber\\\nonumber
&=& \displaystyle\sum^{n}_{i=0} \sum^{n-i}_{j=0}\binom{n + 1}{n-i-j}
(-k)^j(d-1)^i,
\end{eqnarray}
\noindent where in the last equality we have used the fact that $\displaystyle \int_{\CIPT}c_1(\OO_{\CIPT}(1))^{n} = 1$.

By hypothesis, the singularities of $\fol$ are non-degenerates. Then, the number $\# \left [\singf \cap {\CIPT \backslash \ddd}\right ]$ corresponds to the sum of the numbers of Milnor of the singular points of $\fol$ in $\CIPT \backslash \ddd$. Moreover,   for all $p \in \Sing(\fol)\cap \ddd_{reg}$ we have  $Ind_{log\, \ddd,p}(\fol)=0$, since the singularities are non-degenerates.
 Thus,  it follows from Theorem $1$  that 
\begin{eqnarray}\nonumber
\# \left [\singf \cap {\CIPT \backslash \ddd}\right ] 
&=& \displaystyle\sum^{n}_{i=0} \sum^{n-i}_{j=0}\binom{n + 1}{n-i-j}
(-k)^j(d-1)^i.
\end{eqnarray}
\noindent Now, the conclusion of proof can readily be obtained by the signal study of $\delta(k,d,n)$ that was done in Lemma \ref{lemmm}.
$\square$
\end{demostracion}

Particularly, the items (1b) and (2b) of Theorem $3$ characterize the situations in which all the singularities of $\fol$ occur in the  invariant hypersurface $\ddd$. We will present optimal examples.

\begin{exe} Let $\fol$ be the foliation on $\mathbb{P}^3$ induced by the polynomial vector field
\begin{eqnarray}\nonumber
\displaystyle v &=& (-z_1^{k-1} - z_2^{k-1}- z_3^{k-1})\frac{\partial}{\partial z_0} + (z_0^{k-1} - z_2^{k-1}- z_3^{k-1})\frac{\partial}{\partial z_1} +\\\nonumber\\\nonumber &&+ (z_0^{k-1} + z_1^{k-1}- z_3^{k-1})\frac{\partial}{\partial z_2} + (z_0^{k-1} + z_1^{k-1}+ z_2^{k-1})\frac{\partial}{\partial z_3}.
\end{eqnarray}
\noindent The hypersurface $\ddd =\{z_0^{k} +z_1^{k} +z_2^{k} +z_3^{k}=0\}$  is    invariant by $\fol$. 
 It is not difficult see that  $\Sing(\fol) \subset \ddd$.
Note that  $\deg(\ddd) = k$ and $\deg(\fol) = k-1$,  according to item (1.b) of the Theorem 3.
\end{exe}

\begin{exe}
Consider the  foliation $\fol$ induced  by  the vector field $v =  \partial/\partial z_0$.  For each $1\leq i\leq n$, the hypersurface  $\ddd_i = \{z_i = 0\}$ is invariant by $\fol$. Moreover,   that  for all $i=1, \dots, n$, we have
$$\Sing(\fol) = \{(1:0:\ldots:0)\} \subset \ddd_i.$$
Note that we have $\deg(\fol) = 0$ and $\deg(\ddd_i) = 1$, for all $i$. Therefore, if we consider $n$ even, we are  in the case of item (2b) of Theorem 3.

\end{exe}

\end{document}